\documentclass[reqno]{amsart}%
\usepackage{amsfonts}
\usepackage{verbatim}
\usepackage{xcolor}
\usepackage{amsmath}
\usepackage{amssymb}
\usepackage{graphicx}
\usepackage{accents}%
\providecommand{\U}[1]{\protect\rule{.1in}{.1in}}
\newtheorem{theorem}{Theorem}[section]

\numberwithin{equation}{section}
\begin{document}
\title[Hyperbolic systems]{Method for solving hyperbolic systems with initial data on non-characteristic
manifolds with applications to the shallow water wave equations}
\author{Alexei Rybkin}
\address{Department of Mathematics and Statistics, University of Alaska Fairbanks, PO
Box 756660, Fairbanks, AK 99775}
\email{arybkin@alaska.edu}
\thanks{The author is supported in part by the NSF grant DMS-1716975.}
\date{January, 2019}
\subjclass{58J45, 76B15, 35Q35}
\keywords{hyperbolic systems, shallow water wave equations, Carrier-Greenspan transform}

\begin{abstract}
We are concerned with hyperbolic systems of order-one linear PDEs originated
on non-characteristic manifolds. We put forward a simple but effective method
of transforming such initial conditions to standard initial conditions (i.e.
when the solution is specified at an initial moment of time). We then show how
our method applies in fluid mechanics. More specifically, we present a
complete solution to the problem of long waves run-up in inclined bays of
arbitrary shape with nonzero initial velocity.

\end{abstract}
\maketitle

\section{Introduction}

Let $A_{i}=A_{i}\left(  \mathbf{x},t\right)  ,i=0,1,...,n$ be known $m\times
m$ matrix-valued functions of $\mathbf{x}=\left(  x_{i}\right)  \in
\mathbb{R}^{n}$ and $t\in\mathbb{R}$; $\mathbf{f}\left(  \mathbf{x},t\right)
,\mathbf{g}\left(  \mathbf{x}\right)  $ are known $\mathbb{R}^{m}$-valued
functions; and $\mathbf{u}\left(  \mathbf{x},t\right)  $ is an unknown
$\mathbb{R}^{m}$-valued function. We are concerned with the Cauchy problem on
manifolds for hyperbolic systems of linear first-order PDEs%
\begin{align}
\mathbf{u}_{t}+\sum_{i=1}^{n}A_{i}\mathbf{u}_{x_{i}}+A_{0}\mathbf{u} &
=\mathbf{f}\text{ in }\mathbb{R}^{n+1},\label{hypo system}\\
\left.  \mathbf{u}\right\vert _{\Gamma} &  =\mathbf{g},\label{IC}%
\end{align}
where%
\begin{equation}
\Gamma=\{(\mathbf{x},\tau(\mathbf{x}))|\mathbf{x}\in\text{ }\mathbb{R}%
^{n}\}\label{curve}%
\end{equation}
is a n-dimensional manifold in $\mathbb{R}^{n+1}$ defined by a known scalar
function $\tau(\mathbf{x})$. In particular, if $n=1$ then $\Gamma$ is a curve.

Recall that the system (\ref{hypo system}) is called hyperbolic if the matrix%
\begin{equation}
A\left(  \mathbf{x},t,\mathbf{\xi}\right)  :=\sum_{i=1}^{n}\xi_{i}A_{i}\left(
\mathbf{x},t\right)  ,\ \ \mathbf{\xi}=\left(  \xi_{i}\right)  \in
\mathbb{R}^{n},\label{A}%
\end{equation}
is diagonalizable for each $\mathbf{x},\mathbf{\xi}\in\mathbb{R}^{n}$ and
$t\in\mathbb{R}$ with real eigenvalues (playing roles of characteristic
velocities) (see, e.g. \cite{EvansBook}).

If $\tau\left(  \mathbf{x}\right)  =0$ then (\ref{IC}) becomes the standard
initial condition (IC)\ $\left.  \mathbf{u}\right\vert _{t=0}=\mathbf{g}$ for
which the initial value problem (IVP) is well-posed for a large class of
$A_{i},\mathbf{f},\mathbf{g}$ and is very well studied. If $\tau\left(
\mathbf{x}\right)  \neq0$ then (\ref{hypo system})-(\ref{IC}) is well-posed
\cite{EvansBook} only if%
\begin{equation}
\det\left\{  I-A\left(  \mathbf{x},t,\nabla\tau\left(  \mathbf{x}\right)
\right)  \right\}  \neq0\ \ \text{in }\mathbb{R}^{n+1}.\label{det cond}%
\end{equation}
Such manifolds are called non-characteristic, i.e. tangent planes to $\Gamma$
and to characteristic manifolds are nowhere parallel. The problem however is
that solving (\ref{hypo system})-(\ref{IC}) becomes a notoriously hard problem
even for constant $A_{i}$. On the other hand, such problems routinely appear
while linearizing certain quasi-linear PDEs by the hodograph transform, a
technique involving switching independent and dependent variables (see e.g.
\cite{ClarksonetalHodograph89} for an extensive list of PDEs linearizable by
the hodograph and \cite{Qu} for the current literature). As it is simply put
in \cite{Johnson} \textquotedblleft Interchanging the dependent and
independent variables simplifies the governing equations, but complicates the
boundary/initial conditions." This is a serious issue deterring from using the
hodograph transform and when they do use it the problem of dealing with
entangled initial/boundary conditions is solved on an ad hoc basis relying on
various intuitive assumptions with no error analysis.

In the present note we put forward an elementary but very effective method of
converting conditions on manifolds to standard conditions. Once it is done,
one can proceed with any applicable method of solving the standard IBVP. Note
our approach easily extends to higher order linear hyperbolic PDEs since they
can be written as systems. We demonstrate effectiveness of our approach on the
run-up problem of tsunami waves, a well known problem of fluid mechanics, by
improving on (in chronological order) \cite{Yeh03}, \cite{KANOGLU04},
\cite{Kanoglu06}, \cite{Antuono10}, and \cite{Nicolsky18} (see also the
literature cited therein).

We owe the following comment to one of the referees. The book
\cite{Sachdev2000} provides more examples of the applicability of our work to
other nonlinear hyperbolic equations. In particular, to Section 2.2 and
Chapter 5 of \cite{Sachdev2000} where our method could be used to eliminate
some of the simplifications used to form initial conditions post applying a
hodograph transform.

\section{\label{Data projection}The method of data projection}

In this section we introduce what we call the\emph{ }method of data
projection. For simplicity we assume that $A_{i}$ depend only on $\mathbf{x}$
and the nonhomogeneous term $\mathbf{f}$ is absent. It will be clear however
that the general case merely results in more complicated formulas. Introduce
the following differential operator
\[
D\mathbf{u}:=\left[  \sum_{i=1}^{n}A_{i}\left(  \mathbf{x}\right)
\frac{\partial}{\partial x_{i}}+A_{0}\left(  \mathbf{x}\right)  \right]
\mathbf{u.}%
\]
Our IVP then reads%
\begin{equation}
\left\{
\begin{array}
[c]{c}%
\mathbf{u}_{t}=-D\mathbf{u}\\
\mathbf{u}|_{\Gamma}=\mathbf{g}%
\end{array}
\right.  .\label{IVP}%
\end{equation}
Our idea is, given any accuracy $\epsilon$, we find a standard initial
condition $\mathbf{u}|_{t=0}=\widetilde{\mathbf{g}}$ such that the solution
$\widetilde{\mathbf{u}}$ to
\begin{equation}
\left\{
\begin{array}
[c]{c}%
\mathbf{u}_{t}=-D\mathbf{u}\\
\mathbf{u}|_{t=0}=\widetilde{\mathbf{g}}%
\end{array}
\right.  \label{newIVP}%
\end{equation}
is within $O\left(  \varepsilon\right)  $ from the actual solution to
(\ref{IVP}) uniformly in $\left(  \mathbf{x},t\right)  $ in the domain of
interest. I.e., the problems (\ref{IVP}) and (\ref{newIVP}) are equivalent up
to $O\left(  \varepsilon\right)  $.

We call the map $\mathbf{g}\longrightarrow$ $\widetilde{\mathbf{g}}$ the\emph{
}data projection\emph{ }\ of the manifold $\Gamma$ onto the hyperplane
$R^{n}\times\left\{  t=0\right\}  $. Finding $\widetilde{\mathbf{g}}$ is based
upon the Taylor formula in "reverse".

Apply the Taylor formula in one variable $t$ to the (still unknown) solution
$\mathbf{u}\left(  \mathbf{x},t\right)  $ of (\ref{IVP}). For each fixed point
$(\mathbf{x},t)$ we then have
\begin{equation}
\mathbf{u}(\mathbf{x},0)=\sum_{k=0}^{j}\frac{1}{k!}\frac{\partial
^{k}\mathbf{u}(\mathbf{x},t)}{\partial t^{k}}(-t)^{k}+E_{j}\left(
\mathbf{x},t\right)  \label{taylor}%
\end{equation}
with some error $E_{j}$. I.e., (\ref{taylor}) is the Taylor formula about
$\left(  \mathbf{x},t\right)  $ evaluated at $\left(  \mathbf{x},0\right)  $
(not the other way around). Taking in (\ref{taylor}) $t=\tau\left(
\mathbf{x}\right)  $ yields%
\begin{align}
\left.  \mathbf{u}\left(  \mathbf{x},t\right)  \right\vert _{t=0} &
=\sum_{k=0}^{j}\frac{(-\tau(\mathbf{x}))^{k}}{k!}\left.  \frac{\partial
^{k}\mathbf{u}(\mathbf{x},t)}{\partial t^{k}}\right\vert _{\Gamma}+\left.
E_{j}\right\vert _{\Gamma}\nonumber\\
&  =:\mathbf{g}_{j}(\mathbf{x})+\left.  E_{j}\right\vert _{\Gamma
}.\label{nth order projection}%
\end{align}
We call $\mathbf{g}_{j}$ the $j^{\text{th}}$ order projection of initial data
$\mathbf{g}$ onto $R^{n}\times\left\{  t=0\right\}  $. One can see now that if
we are able to properly compute all $\left.  \dfrac{\partial^{k}%
\mathbf{u}(\mathbf{x},t)}{\partial t^{k}}\right\vert _{\Gamma}$ then
$\mathbf{g}_{j}$ produces a desirable standard IC $\widetilde{\mathbf{g}}$ .
Indeed, given error $\varepsilon$ (no matter how small), we take $j$ so large
as $\left\vert E_{j}\left(  \mathbf{x},t\right)  \right\vert <\varepsilon$
uniformly in the domain of interest and hence%
\[
\mathbf{u}|_{t=0}=\mathbf{g}_{j}+O\left(  \varepsilon\right)  .
\]
Thus the solution $\widetilde{\mathbf{u}}$ to (\ref{newIVP}) with
$\widetilde{\mathbf{g}}=\mathbf{g}_{j}$ will coincide with the solution
$\mathbf{u}$ of the original problem (\ref{IVP}) up to $O\left(
\varepsilon\right)  $.

We now show how to compute the Taylor coefficients in
(\ref{nth order projection}). The zeroth one is obvious%
\[
\mathbf{g}_{0}(\mathbf{x})=\mathbf{u}(\mathbf{x},t)|_{\Gamma}=\mathbf{g}%
\left(  \mathbf{x}\right)  \text{ \ (IC in (\ref{IC})).}%
\]
All of the other Taylor coefficients in (\ref{nth order projection}) can also
be explicitly computed. Start with the first-order one:
\begin{equation}
\left.  \mathbf{u}_{t}\right\vert _{\Gamma}=-\left.  (D\mathbf{u})\right\vert
_{\Gamma}.\label{1st coeff}%
\end{equation}
We now compute $\left.  (D\mathbf{u})\right\vert _{\Gamma}$ by the chain rule.
One has%
\begin{align*}
D\mathbf{g} &  =\left[  \sum_{i=1}^{n}A_{i}\frac{\partial}{\partial x_{i}%
}+A_{0}\right]  \mathbf{g}\\
&  =\sum_{i=1}^{n}A_{i}\left(  \left.  \mathbf{u}_{x_{i}}\right\vert _{\Gamma
}+\tau_{x_{i}}\left.  \mathbf{u}_{t}\right\vert _{\Gamma}\right)
+A_{0}\mathbf{g}\qquad\text{(by the chain rule)}\\
&  =\left.  \left(  \sum_{i=1}^{n}A_{i}\mathbf{u}_{x_{i}}+A_{0}\mathbf{u}%
\right)  \right\vert _{\Gamma}+\sum_{i=1}^{n}\tau_{x_{i}}A_{i}\left.
\mathbf{u}_{t}\right\vert _{\Gamma}\\
&  =\left.  \left(  D\mathbf{u}\right)  \right\vert _{\Gamma}-\left(
\sum_{i=1}^{n}\tau_{x_{i}}A_{i}\right)  \left.  \left(  D\mathbf{u}\right)
\right\vert _{\Gamma}\text{\ \ \ \ \ (by (\ref{1st coeff}))}\\
&  =\left(  I-\sum_{i=1}^{n}\tau_{x_{i}}A_{i}\right)  \left.  \left(
D\mathbf{u}\right)  \right\vert _{\Gamma}.
\end{align*}
It follows that%
\[
\left.  \left(  D\mathbf{u}\right)  \right\vert _{\Gamma}=\left(  I-\sum
_{i=1}^{n}\tau_{x_{i}}A_{i}\right)  ^{-1}D\mathbf{g},
\]
and thus for the first order Taylor coefficient we finally have
\begin{equation}
\left.  \mathbf{u}_{t}\right\vert _{\Gamma}=-\left(  I-A\right)  ^{-1}D\left(
\mathbf{u}|_{\Gamma}\right)  =-\left(  I-A\right)  ^{-1}D\mathbf{g}%
,\label{partial of t}%
\end{equation}
where, recalling (\ref{A}),%
\[
A:=\sum_{i=1}^{n}\tau_{x_{i}}\left(  \mathbf{x}\right)  A_{i}\left(
\mathbf{x}\right)  =A\left(  \mathbf{x},\nabla\tau\left(  \mathbf{x}\right)
\right)  .
\]

Our computation of higher order Taylor coefficients will be based on the
following observation. Since $\partial/\partial t$ and $D$ commute,
$\mathbf{u}_{t}$ is also a solution to
\[
\left(  \mathbf{u}_{t}\right)  _{t}=-D\left(  \mathbf{u}_{t}\right)
\]
with the initial condition given by (\ref{partial of t})%
\[
\left.  \mathbf{u}_{t}\right\vert _{\Gamma}=-\left(  I-A\right)
^{-1}D\mathbf{g}.
\]
Thus, if $\mathbf{u}$ is the solution originated from $\mathbf{g}$ then
$\mathbf{u}_{t}$ is the solution originated from the initial data $-\left(
I-A\right)  ^{-1}D\mathbf{g}$. By induction one concludes that $\dfrac
{\partial^{k}\mathbf{u}}{\partial t^{k}}$ is the solution originated from
$\left[  -\left(  I-A\right)  ^{-1}D\right]  ^{k}\mathbf{g}$, $k=2,3,4,...$%
\newline Therefore, we get the following nice formula
\[
\left(  \left.  \dfrac{\partial^{k}\mathbf{u}}{\partial t^{k}}\right\vert
_{\Gamma}\right)  =\left(  -\left(  I-A\right)  ^{-1}D\right)  ^{k}%
\mathbf{g},\ \ k=0,1,2,\ ...
\]
Substituting this into (\ref{nth order projection}) one has%

\begin{equation}
\mathbf{g}_{j}=\mathbf{g}+\sum_{k=0}^{j}\frac{\tau^{k}}{k!}\left(  \left(
I-A\right)  ^{-1}D\right)  ^{k}\mathbf{g} \label{final nth order projection}%
\end{equation}
and we finally arrive at

\begin{theorem}
Let $\Gamma=\{(\mathbf{x},\tau(\mathbf{x}))|\mathbf{x}\in$ $\mathbb{R}^{n}\}$
be a non-characteristic manifold for the hyperbolic system%
\begin{equation}
\mathbf{u}_{t}+\sum_{i=1}^{n}A_{i}\left(  \mathbf{x}\right)  \mathbf{u}%
_{x_{i}}+A_{0}\left(  \mathbf{x}\right)  \mathbf{u}=0\text{ in }%
\mathbb{R}^{n+1},\label{system}%
\end{equation}
and%
\begin{equation}
\left.  \mathbf{u}\right\vert _{\Gamma}=\mathbf{g}\left(  \mathbf{x}\right)
.\label{ic}%
\end{equation}
Fix accuracy $\varepsilon$ and let $j$ be chosen to satisfy ($\left\Vert
\cdot\right\Vert $ is the Euclidean norm in $\mathbb{R}^{m}$)
\[
\max_{\mathbf{x}\in\mathbb{R}^{n}}\left\Vert \frac{\tau^{j+1}\left(
\mathbf{x}\right)  }{\left(  j+1\right)  !}\left\{  \left[  I-A\left(
\mathbf{x},\nabla\tau\left(  \mathbf{x}\right)  \right)  \right]
^{-1}D\right\}  ^{j+1}g\left(  \mathbf{x}\right)  \right\Vert <\varepsilon,
\]%
\[
A\left(  \mathbf{x},\nabla\tau\left(  \mathbf{x}\right)  \right)  :=\sum
_{i=1}^{n}\tau_{x_{i}}\left(  \mathbf{x}\right)  A_{i}\left(  \mathbf{x}%
\right)  ,\ \ \ D:=\sum_{i=1}^{n}A_{i}\left(  \mathbf{x}\right)
\frac{\partial}{\partial x_{i}}+A_{0}\left(  \mathbf{x}\right)  .
\]
If $\widetilde{\mathbf{u}}$ is the solution to (\ref{system}) subject to the
initial condition%
\[
\left.  \widetilde{\mathbf{u}}\right\vert _{t=0}=\mathbf{g}+\sum_{k=0}%
^{j}\frac{\tau^{k}}{k!}\left(  \left(  I-A\right)  ^{-1}D\right)
^{k}\mathbf{g}%
\]
then the solution $\mathbf{u}$ to (\ref{system})-(\ref{ic}) is subject to
$\mathbf{u}\left(  \mathbf{x},t\right)  =\widetilde{\mathbf{u}}\left(
\mathbf{x},t\right)  +O\left(  \varepsilon\right)  .$
\end{theorem}

The map $\mathbf{g}\longrightarrow\mathbf{g}_{j}$ is linear and well-defined
as long as the matrix $I-A$ is non-singular, which holds if (\ref{det cond})
does. An important feature of our data projection method is that, by
construction, both $\mathbf{u}$ and $\widetilde{\mathbf{u}}$ solve (exactly)
the same system (\ref{hypo system}) but satisfy different (equivalent) IC
conditions $\mathbf{g}$ and $\widetilde{\mathbf{g}}=\mathbf{g}_{j}$. Of course
$\mathbf{u}$ and $\widetilde{\mathbf{u}}$ can be made as close as one wishes
(while $\mathbf{g}$ and $\widetilde{\mathbf{g}}$ need not be close at all). \ 

\section{Applications to the non-linear shallow water wave system}

In this section we apply our formalism to the hyperbolic 1+1 quasi-linear
shallow water system (see, e.g. \cite{Lannes13}) modeling e.g. the tsunami
wave run-up and run-down. For the so-called inclined bathymetries this system
reads (in dimensionless variables)%
\begin{equation}
\left\{
\begin{array}
[c]{ccc}%
\eta_{t}+\left(  1+\eta_{x}\right)  u+c\left(  x+\eta\right)  u_{x}=0 &  &
\text{(continuity equation)}\\
u_{t}+uu_{x}+\eta_{x}=0 &  & \text{(momentum equation)}\\
\eta\left(  x,0\right)  =\eta_{0}\left(  x\right)  ,\ \ u\left(  x,0\right)
=u_{0}\left(  x\right)   &  & \text{(initial conditions)}%
\end{array}
\right.  ,\label{SWE1}%
\end{equation}
where $\eta$ is the unknown water elevation over unperturbed level, $u$ is the
unknown cross-section averaged flow velocity, and $c\geq0$ is a known function
solely encoding the information about the shape of our inclined bathymetry
\cite{Rybkin}, \cite{Raz2017}. The point $x=0$ in (\ref{SWE1}) corresponds to
the unperturbed coast line and the $x$-axis is directed off shore. The main
feature of this problem is a moving (wet/dry) boundary also know as
run-up/run-down. The substitution%
\begin{equation}%
\begin{array}
[c]{cc}%
\varphi\left(  \sigma,\tau\right)  =u\left(  x,t\right)  , & \psi\left(
\sigma,\tau\right)  =\eta\left(  x,t\right)  +u^{2}\left(  x,t\right)  /2,\\
\sigma=x+\eta\left(  x,t\right)  , & \tau=t-u\left(  x,t\right)  ,
\end{array}
\label{CG}%
\end{equation}
referred to as the Carrier-Greenspan (CG)\ transform or CG hodograph, turns
(\ref{SWE1}) into the linear (strictly) hyperbolic system%
\begin{align}
&  \left\{
\begin{array}
[c]{c}%
\mathbf{\phi}_{\tau}+A(\sigma)\mathbf{\phi}_{\sigma}+B\mathbf{\phi}=0\\
\left.  \mathbf{\phi}\right\vert _{\Gamma}=\mathbf{\phi}_{0}\left(
\sigma\right)
\end{array}
\right.  ,\label{SWE in matrix form}\\
A\left(  \sigma\right)   &  =\left(
\begin{array}
[c]{cc}%
0 & 1\\
c\left(  \sigma\right)   & 0
\end{array}
\right)  ,\ \ B=\left(
\begin{array}
[c]{cc}%
0 & 0\\
1 & 0
\end{array}
\right)  ,\ \ \mathbf{\phi}=\left(
\begin{array}
[c]{c}%
\varphi\\
\psi
\end{array}
\right)  .\label{A, B, etc}%
\end{align}
where $\Gamma$ and $\mathbf{\phi}_{0}$ will be given later. The CG transform
was introduced in \cite{Carrier1} and has become a standard tool in the study
of the run-up/run-down process (see, e.g. \cite{Didenkulova11b},
\cite{Rybkin}, \cite{Synolakis87} and the extensive literature cited therein).
The form (\ref{CG}) is taken from our \cite{Raz2017}. Besides linearizing
(\ref{SWE1}), the CG turns the moving boundary into the fixed point $\sigma
=0$. It has however a serious drawback: the IC in (\ref{SWE in matrix form})
is no longer standard. Indeed, under (\ref{CG}), the (horizontal) line $t=0$
in the plane $\left(  x,t\right)  $ becomes the parametric curve
$\Gamma=\left(  x+\eta_{0}(x),-u_{0}(x)\right)  $ in the $(\sigma,\tau)$ plane
and, denoting the inverse of $x+\eta_{0}(x)$ by $\gamma\left(  \sigma\right)
$, one has
\begin{align}
\Gamma &  =\left\{  \left(  \sigma,-u_{0}\left(  \gamma\left(  \sigma\right)
\right)  \right)  |\ \sigma\geq0\right\}  ,\label{Gamma}\\
\left.  \left(
\begin{array}
[c]{c}%
\varphi\\
\psi
\end{array}
\right)  \right\vert _{\Gamma} &  =\left.  \left(
\begin{array}
[c]{c}%
u_{0}\\
\eta_{0}+u_{0}^{2}/2
\end{array}
\right)  \right\vert _{\gamma\left(  \sigma\right)  }=:\left(
\begin{array}
[c]{c}%
\varphi_{0}\left(  \sigma\right)  \\
\psi_{0}\left(  \sigma\right)
\end{array}
\right)  =:\mathbf{\phi}_{0},\label{IC transformed}%
\end{align}
which defines $\mathbf{\phi}_{0}$ in (\ref{SWE in matrix form}). We see from
(\ref{Gamma}) that the IC (\ref{IC transformed}) is standard iff the initial
velocity $u_{0}=0.$ The latter has been a typical assumption in much of the
previous literature (which however does not give a valid inundation picture
caused by a tsunami wave). That is why it has been a good open problem since
\cite{Carrier1} how to make the CG transform run for general IC. We refer to
\cite{Yeh03}, \cite{KANOGLU04}, \cite{Kanoglu06}, \cite{Antuono10} where this
problem was addressed under certain assumptions of relative smallness of the
IC (and only in the context of the plane beach). We discovered our method in
\cite{Nicolsky18} in a particular case while trying to put previous works on a
solid footing.

We now apply the method of data projection to (\ref{SWE in matrix form}).
Given accuracy $\varepsilon$, choose $j$ so that%
\[
\max_{\sigma\geq0}\left\Vert \frac{\varphi_{0}^{j+1}\left(  \sigma\right)
}{\left(  j+1\right)  !}\left\{  \left[  I+\varphi_{0}^{\prime}\left(
\sigma\right)  A\left(  \sigma\right)  \right]  ^{-1}\left[  A\left(
\sigma\right)  \frac{d}{d\sigma}+B\right]  \right\}  ^{j+1}\mathbf{\phi}%
_{0}\left(  \sigma\right)  \right\Vert <\varepsilon.
\]
and construct the $j$th order projection%
\begin{equation}
\mathbf{\phi}_{j}\left(  \sigma\right)  =\mathbf{\phi}_{0}\left(
\sigma\right)  +\sum_{k=1}^{j}\frac{\left(  -\varphi_{0}\left(  \sigma\right)
\right)  ^{k}}{k!}\left\{  \left[  I+\varphi_{0}^{\prime}\left(
\sigma\right)  A\left(  \sigma\right)  \right]  ^{-1}\left[  A\left(
\sigma\right)  \frac{d}{d\sigma}+B\right]  \right\}  ^{k}\mathbf{\phi}%
_{0}\left(  \sigma\right)  . \label{Phi sub n}%
\end{equation}
We can then solve (\ref{SWE in matrix form}) with the standard IC $\left.
\mathbf{\phi}\right\vert _{\tau=0}=\mathbf{\phi}_{j}\left(  \sigma\right)  $
by any suitable method. In fact, for the so-called power shaped bays,
$c\left(  \sigma\right)  \sim\sigma$ and it can be solved by the Hankel
transform in terms of Bessel functions \cite{Nicolsky18}. Performing the
inverse CG transform (\ref{CG}) solves the original problem (\ref{SWE1}) in
the physical space. The latter is, in general, not explicit but can easily be
done numerically without affecting the total accuracy, which remains $O\left(
\varepsilon\right)  .$ In fact, we can call our method exact as the error it
introduces can be made negligible comparing with the one inherited by the
shallow water approximation leading to the very system (\ref{SWE1}).

Note that while (\ref{Phi sub n}) looks unwieldy (the reader is invited to
amuse him/herself with trying to unzip it even for $j=1$), its numerical
implementation is not a problem. It was the matrix form
(\ref{SWE in matrix form}) that made our derivation quite transparent.

We emphasize that an important feature of our method is that the procedure
works smoothly as long as the wave non-breaking condition \cite{Rybkin} is
satisfied and no extra assumptions of smallness (made in the previous
literature)\ are needed. Extensive numerical verification and simulations
(which will be published elsewhere) show that our method is robust and can be
effectively used for rapid forecasting of characteristics of inundation zone
due to large-amplitude sea waves.

In the conclusion, we mention that our method can be adapted to treat boundary
conditions and improve on the relevant results of \cite{Anderson17},
\cite{Antuono10}, \cite{Harrisetal16}, \cite{Synolakis87}. We plan to address
this elsewhere.

\section{Acknowledgments}

We are grateful to Efim Pelinovsky for stimulating email discussions and an
annonemous referee for useful literature hints.

\end{document}